\input amstex


\def\b1{\text{\bf 1}}

\def\CC{{\Cal C}}
\def\CD{{\Cal D}}

\def\CE{{\Cal E}}

\def\CI{{\Cal I}}

\def\CM{{\Cal M}}

\def\CO{{\Cal O}}
\def\CP{{\Cal P}}

\def\CT{{\Cal T}}

\def\gr{\text{gr}}

\def\Ext{\text{Ext}}

\def\#{\,\check{}}

\def\Coker{\text{Coker}}
\def\Ker{\text{Ker}}


\def\imbed{\hookrightarrow}

\def\hra{\hookrightarrow}
\def\iso{\buildrel\sim\over\rightarrow}

\parskip=6pt

\def\A{{\Cal A}}
\def\B{{\Cal B}}
\def\D{{\Cal D}}
\def\I{{\Cal I}}
\def\T{{\Cal T}}
\def\P{{\Cal P}}
\def\O{{\Cal O}}

\def\mapsto{{\rightarrow}}

\def\Pone{{\Bbb P}^1}

\documentstyle{amsppt}
\document
\magnification=1100
\NoBlackBoxes

\centerline {TILTING EXERCISES}
\medskip
\centerline {A.~Beilinson,
R.~Bezrukavnikov, and I.~Mirkovi\'c\footnote{All authors are partially 
supported by NSF grants:
A.B. is supported by grant DMS-0100108, and
R.B. by grant DMS-0071967.}}

\bigskip

 This is a geometry-oriented review of
the basic formalism of tilting objects (originally due to Ringel, 
see [Ri]\footnote{Notice that our terminology differs from the one in
{\it loc. cit.} and some other sources; the term ``tilting'' there is used for 
a weaker property.}, \S 5).
 In the
first section we explain that tilting
extensions form a natural framework for the
gluing construction from [B1] and [MV]. We
show that in case of a stratification with
contractible strata, the homotopy category of
complexes of tilting perverse sheaves is
equivalent to the derived category of sheaves
smooth along the stratification. Thus tilting
objects play the role similar to projective or
injective ones (with advantage of being
self-dual and having local origin).  In the
second section we discuss tilting perverse
sheaves smooth along the Schubert
stratification of the flag space (or,
equivalently, tilting objects in the
Bernstein-Gelfand-Gelfand category
$\O$). In this case  a Radon transform interchanges
tilting, projective, and injective modules. 
 As a corollary, we give a short
proof of Soergel's Struktursatz [S1], and
describe the Serre functor for
$D^b (\O )$ (as conjectured by M.~Kapranov).

 We refer to
[M]  for a much more thorough exposition of many
other aspects of the theory.

\medskip

This article is a modest present to Borya
Feigin -- with love, and sadness to see him
so rarely these days. 

{\it ..Mr.~Fagin took the opportunity of
reading Oliver a long lecture on the crying 
sin of ingratitude: of which he clearly
demonstrated he had been guilty,  to no
ordinary extent, in wilfully absenting himself
from the society of his  anxious friends..}

\bigskip

\centerline{\bf \S 1 Generalities.}

\medskip

We consider algebraic varieties over an
algebraically closed field $k$. Below
``perverse sheaf" means either plain perverse
$\Bbb Q_l$-sheaf, $l\ne char(k)$, or perverse
sheaf with respect to classical
topology with coefficients in any field of
characteristic 0 (if $k=\Bbb C$), or  holonomic
$\CD$-module (in case $char (k)=0$). For a
variety $X$ we denote by
$\CM (X)$ the abelian category of perverse
sheaves on
$X$, and by $D (X)$ its bounded
derived category (which is the same as the
usual ``topological" derived category of
constructible
$\Bbb Q_\ell$-complexes or complexes of
$\CD$-modules with holonomic cohomology, see
[B2]).

\medskip

{\bf 1.1} Let $X$ be an algebraic
variety,
$i:Y\hra X$ a closed subvariety,
$j:U:=X\smallsetminus Y\hra X$ the
complementary embedding. Let $M$ be a perverse
sheaf on
$X$; then we have\footnote{Here the
t-structure is of the middle perversity.}
$i^! M\in D(Y)^{\ge 0}$,
$i^* M \in D(Y)^{\le 0}$. We
say that $M$ is a {\it tilting} perverse sheaf with
respect to
$Y$ (or a {\it
$Y$-tilting} perverse sheaf) if both
$i^! M$,
$i^* M$ are perverse sheaves. The standard exact triangles
together with left (respectively, right) exactness of $j_*$, $j_!$
show that $M$ is tilting if and only if both
$j_*j^* M$ and $j_! j^* M$ are perverse sheaves and the
canonical morphisms $M\to j_* j^* M$, $j_! j^* M \to M$
are, respectively, surjective and injective. The category
of $Y$-tilting perverse sheaves is closed under extensions
and  Verdier duality.

{\bf Proposition.} Let $M_U$ be a perverse sheaf on $U$
such that both $j_* (M_U )$ and $j_! (M_U )$
are perverse sheaves on $X$. Then there exists
a
$Y$-tilting perverse sheaf
$M$ on
$X$ such that
$M|_U =M_U$. We call such $M$ a {\it 
$Y$-tilting extension} of $M_U$ to $X$.

{\it Proof.} (a) Set $A:= \Ker (j_! ( M_U )\to
j_* ( M_U ))$,
$B:=
\Coker (j_! ( M_U )\to j_* (M_U ) )$. These are
perverse sheaves supported on $Y$. Let
$c\in\Ext^2 (B,A)$ be the Yoneda class of
exact sequence
$0\to A\to j_! ( M_U )\to j_* (M_U )\to B\to
0$.

(b) If $c$
vanishes, then there exists a
perverse sheaf $M$ together with a 3 step
filtration $M_0 \subset M_1 \subset M$
and identifications
$M_0 =A$, $M_1 =j_! (M_U )$, $M/M_0
=j_* (M_U )$,
$M/M_1 =B$ compatible in the obvious manner with the above
exact sequence. Then from the exact sequences 
$$0\to A\to M\to j_*(M_U)\to 0,$$
$$0\to j_!(M_U)\to M \to B\to 0$$
we see that $A\iso i^!(M)$, $i^*(M)\iso B$; thus 
 $M$ is  $Y$-tilting 
(we call 
it a {\it 
minimal tilting extension} of $M_U$ for obvious reasons),
and we are done.

(c) If $c\neq 0$ then we have to correct our exact
sequence. To do this notice that by [B2] the Yoneda
Ext's are the same as Ext's in the usual derived category
of sheaves on
$X$. The latter can be computed inside the derived
category of sheaves on $Y$, and then as  Yoneda Ext of
perverse sheaves on $Y$. Thus one can find an exact
sequence $0\to A\to C\to D\to B\to 0$ of perverse sheaves
supported on
$Y$ of the  class
$-c$. Let $0
\to A
\to C'
\to D'
\to B\to 0$ be the Baer sum of the two Yoneda extensions.
Its class vanishes, thus there exists a perverse sheaf $M$
together with a 3 step filtration $M_0 \subset M_1
\subset M$ such that $M_0 =A$, $M_1 =C'$, $M/M_0 =D'$,
$M/M_1 =B$ compatible
 in the obvious manner with the above
exact sequence. Since $j_! (M_U )\hra C'\hra M$ we see that $i^*(M)=
\Coker (j_!(M_U)\to M)$ is a perverse sheaf,
and since 
$M\twoheadrightarrow D'\twoheadrightarrow j_* (M_U )$ we see that
$i^!(M)=\Ker(M \twoheadrightarrow j_* (M_U ))$ is perverse; thus
 $M$ is
tilting, and we are done.
\hfill$\square$

\medskip
{\it Remarks.} (i) The conditions of Proposition are always
satisfied if $j$ is an affine embedding.

(ii) If $Y$ is a divisor given by equation $f=0$ then
the ``maximal extension"
$\Xi_f (M_U )$ from [B1] is a functorial tilting
extension.

\medskip
{\bf 1.2} Let us show that the gluing
construction from [B1] and [MV] fits
naturally into
 the setting of tilting extensions.

Let $M_U$ be a perverse sheaf on $U$, and let $M^{tilt}$ be a 
$Y$-tilting extension of $M_U$. Set $\Psi :=i^! M^{tilt}$, $\Psi' :=
i^*M^{tilt}$; let $\tau :\Psi \to
\Psi'$ be the composition of the canonical morphisms $
\Psi \hra M^{tilt}\twoheadrightarrow \Psi'$.

{\bf Proposition.} The category $\CE $ of extensions of
$M_U$ to
$X$ is canonically equivalent to the category $\CC$ of
diagrams
$\Psi
\buildrel\alpha\over\to \Phi
\buildrel\beta\over\to \Psi'$ where $\Phi$ is a perverse
sheaf on $Y$ and morphisms $\alpha$, $\beta$ are such that
$\beta\alpha =\tau$.

{\it Proof.}\footnote{This is an immediate
generalization of the proof of [B1] 3.1 that
dealt with the particular case
$M^{tilt} =\Xi_f (M_U )$.} (a)
The functor $\CE
\to
\CC$ sends $M\in\CE$ to $(\Phi ,\alpha,\beta )\in\CC$
defined as follows. Consider a short complex
$$F=F(M):= (j_! (M_U )\to M\oplus M^{tilt} \to
j_* ( M_U )) \tag 1.2.1$$ where the
differentials are defined by the property that
their restriction to $U$ are, respectively, the
diagonal embedding and the anti-diagonal
projection. Set $\Phi =\Phi (M):= H^0 F$.
Notice that the tilting property of
$M^{tilt}$ assures
$H^{\neq 0} F=0$.
We have the obvious morphisms 
$$\CC one (M^{tilt}
\to j_* (M_U )) [-1]\to F\to \CC one (j_! (M_U
)
\to M^{tilt}). \tag 1.2.2$$ Passing to
cohomology, we get
$\Psi
\buildrel\alpha\over\to \Phi
\buildrel\beta\over\to \Psi'$. It is clear that $(\Phi
,\alpha ,\beta )\in\CC$.

(b) The inverse functor $\CC \to\CE$ sends
$\Phi =(\Phi ,\alpha,\beta )\in\CC$ to
$M=M(\Phi )\in\CE$  defined as follows.
Consider a short complex
$$G=G(\Phi ):=(\Psi
\to
\Phi \oplus M^{tilt} \to
\Psi') \tag 1.2.3$$ where the differentials
are, respectively, the sum of $\alpha$ and the
canonical embedding $\Psi \to M^{tilt}$ and
the difference of $\beta$ and the canonical
projection $M^{tilt}\to\Psi'$. Set $M:= H^0
G$. Notice that $H^{\neq 0} G =0$.

It remains to show that functors from (a) and
(b) are mutually inverse. To identify $M(\Phi
(M))$ with $M$ let us replace $\Psi \to \Phi
\to \Psi'$ in the definition of complex $G(\Phi
)$ by $(1.2.2)$. We get a complex whose
cohomology equals $M(\Phi (M))$. On the other
hand, by construction, this complex carries a
3 step filtration with successive quotients
equal to the cone of the identity morphism of
$\CC one (M^{tilt}\to j_* (M_U )) [-1]$, $M$,
and the cone of the identity morphism of
$\CC one (j_! ( M_U )\to M^{tilt})[-1]$. Thus
its cohomology equals
$M$. The construction of the isomorphism $\Phi
(M(\Phi ))\iso
\Phi$ is similar and left to the
reader.
\hfill$\square$

\medskip

{\it Remark.} It follows from the part (b) of the proof
that $i^* M =\CC one (\alpha : \Psi \to \Phi )$, $i^! M=
\CC one(\beta :
\Phi
\to
\Psi')[-1]$. Thus $M$ is tilting if and only if $\alpha$
is injective and $\beta$ is surjective.

\medskip
{\bf 1.3} Suppose  that our variety $X$ 
carries
a stratification $\{ X_\nu \}$; let
$i_\nu :X_\nu \hra X$ be the locally
closed embeddings of the strata. We say that a
perverse sheaf
$M$ is {\it tilting} with respect to our
stratification if for every $\nu$ both
complexes $i_\nu^! M, i^*_\nu M$ are perverse
sheaves on $X_\nu$.

Assume that each $i_\nu$ is an affine
embedding.

{\bf Proposition.}
A perverse sheaf $M$
is tilting with respect to our stratification if and only
if it satisfies the following two
conditions:

1. $M$ can be
represented as a successive extension of perverse sheaves
of type $i_{\nu *}N_\nu$ where $N_\nu$ is a
perverse sheaf on
$X_\nu$.

2. Same with $i_{\nu *}$ replaced by $i_{\nu !}$.

{\it Proof.} Our conditions obviously imply that $M$ is
tilting (notice that $i^!_\mu i_{\nu *} N_\nu$ equals
$N_\nu$ if $\mu =\nu$ and 0 otherwise). 
Conversely, suppose
that
$M$ is tilting. Choose a closed filtration $X\supset X_1
\supset ..\supset X_n \supset X_{n+1}=\emptyset$ such that
$X_i
\smallsetminus X_{i+1}$ is a single stratum. Set $j:
U:= X\smallsetminus X_n \hra X$. Using induction by $n$ we
can assume that $j^* M$ is a successive extension of
perverse sheaves $j^* i_{\nu *}N_\nu $. Thus $j_* j^* M$
is a successive extension of $j_* j^* i_{\nu *}N_\nu
=i_{\nu *}N_\nu$, and the tilting property assures that
$M$ is an extension of $j_* j^* M$ by a perverse sheaf
$i_{n*} i_n^! M$. So condition 1 holds. Condition 2 is
checked in the dual manner.  \hfill$\square$

\medskip
{\bf 1.4} We are in situation of 
 1.3; assume in addition that every
$X_\nu$ is smooth and connected. Let $D=D
(X,\{ X_\nu \} )\subset D (X)$ be the full
subcategory of
 complexes constant along $\{ X_\nu \}$,
i.e., those  $F\in D (X)$ that for every
$\nu$ the complex $i^*_\nu F$ has constant
cohomology sheaves. To assure that $D$ is a
reasonable object to deal with, we assume the
following two properties:

- The cohomology
groups with constant coefficients $H^1 (X_\nu
)$ vanish. Then $D$ is a triangulated
subcategory of $D(X)$. Notice that $D$ is
generated by objects $i_{\nu !} M_\nu$ where
$M_\nu$ are constant (perverse) sheaves on
$X_\nu$.

- One has $i_{\nu *}M_\nu \in D$, i.e.
 $D$ is preserved by the Verdier
duality. Then $D$ is a t-category with
core $\CM =\CM (X,\{ X_\nu \}):= \CM \cap D$;
its irreducible objects are middle
extensions of constant perverse sheaves of
rank 1 on strata. 

Suppose, in addition, that $H^2 (X_\nu
)=0$ for every $\nu$. 

\medskip

{\it Remark.}\footnote{We thank the referee to whom this remark is due.}
 Under the above assumptions the category $\CM$
is what different authors call an abstract Kazhdan-Lusztig category,
or a highest weight category, or a quasi-hereditary category (see e.g.
[BGS], \S 3.2 and reference therein).
 Statements parallel to the next two
Propositions are true  (and apparently well-known to the experts)
for a general category of this sort. 

\medskip

Let
$\CT =\CT (X,\{ X_\nu \}) \subset \CM$ be the
full subcategory of tilting sheaves
with respect to our stratification.

{\bf Proposition.} The support of an
indecomposable object
 $M\in\CT$ is irreducible, i.e., it is the
closure of some stratum $X_{\nu } $, and  $i_\nu^! M $ is
a constant (perverse) sheaf of rank 1 on $X_\nu$. The map
$M\mapsto Supp \, M$ is a bijection between the set of
isomorphism classes of indecomposable objects in
$\CT$ and the set of strata.

{\it Proof.}  Use induction by the
number $n$ of strata. We follow notation of the proof of
Proposition 1.3. By induction our statement is
true for the category $\CT_U$ of tilting
sheaves on
$U$ equipped with the induced stratification. For
every object
$M_U
\in\CT_U$ the complexes $j_* (M_U )$, $j_!
(M_U )$ are perverse sheaves (use 
1.3). The class
$c$ from
 part (a) of the proof of Proposition 1.1
vanishes since $H^2 (X_n )=0$, so $M_U$ admits
a minimal tilting extension
$M\in \CT$ (see ibid., part (b)).  Remark
in 1.2 implies that for
indecomposable
$M_U$ the above $M$ is indecomposable, and
every indecomposable tilting extension of $M_U$
is isomorphic to $M$. It also implies that
every tilting extension of a decomposable
$M_U$ is decomposable. We are done. 
\hfill$\square$

\medskip

{\bf 1.5} We are in situation 1.4, and assume,
in addition, that $H^{>0}(X_\nu
)=0$.

{\bf Proposition.} One has canonical
equivalences of triangulated categories
$$K^b (\CT )\iso D^b \CM \iso D.  \tag 1.5.1$$
Here
$K^b \CT$ is the homotopy category of bounded
complexes in $\CT$.

{\it Proof.} The functors $K^b (\CT )\to D^b
\CM \to D$ in $(1.5.1)$ are the obvious ones. 

 (i) Let us show that the
composition $K^b \CT \to D$ is an equivalence
of categories. By Proposition 1.4
 the image of $K^b \CT$ generates $D$, so
it suffices to prove that 
for every $M,N \in \CT$ one has
$Ext^{>0}_{D}(M,N)=0$. By Proposition 1.3, one
needs to check that
$Ext^{>0}_{D}(M,N)=0$ for $M=i_{\mu !}M_\mu$,
$N=i_{\nu *}N_\nu$, where $M_\mu$,
$N_\nu$ are constant
perverse sheaves  on strata $X_\mu$, $X_\nu$ 
respectively. This follows by adjunction if
$\mu\neq\nu$, and by the vanishing of the
higher cohomology of strata if $\mu =\nu$.

(ii) Let us show  that $D^b \CM \to
D$ is an equivalence of categories. This is a
t-exact functor which identifies the cores, so
it suffices to check that the morphism of the
$\delta$-bifunctors
$Ext^\cdot_{D^b\CM} \to Ext^\cdot_D $ on 
$\CM^\circ \times\CM$ is an isomorphism, or,
equivalently, that $ Ext^\cdot_D $ is
effaceable. 
By (i) our functor $D^b \CM \to
D$ admits a right inverse, so $Ext_D^\cdot$
is a quotient functor of 
$ Ext^\cdot_{D^b \CM}$,  hence it is
effaceable, q.e.d. 
 \hfill$\square$

\medskip

{\it Remark.} An alternative proof of the second equivalence
in (1.5.1) can be found in [BGS], Corollary 3.3.2 
on page 500.  

\bigskip

\centerline{\bf \S 2 The case of Schubert stratification.}

\medskip

{\bf 2.1.} Let $G$ be a semisimple algebraic
group. Let $X=G/B$ be the flag variety
stratified by the Schubert cells $ X_w$, $w\in
W$, where
$W$ is the Weyl group. Our stratified space
satisfies conditions of 1.5. 
Set $ D:=D
(X,\{X_w\} )$, and let $\O\subset D$  be the
category of perverse sheaves.

For $w\in W$  let $L_w, T_w\in \O$ be,
respectively, irreducible and
 indecomposable tilting objects supported
on the closure  of $X_w$;  let
$I_w$ and $P_w$ be, respectively, an injective
hull and projective cover of
$L_w$. Let
$\T$, $\P$,
$\I$ be the categories of, respectively,
tilting, projective, and injective objects. We
also let
$\Delta_w=i_{w!}(M_w)$,
$\nabla_w =i_{w*}(M_w)$ where $M_w$ is the
 constant perverse sheaf of rank 1 on $X_w$.

Let $\O_{> 0}\subset \O$ be  the Serre subcategory generated
by $L_w$, $w\ne e$ (where $e\in W$ is the identity); $\O_0=\O/\O_{>0}$,
 and $\pi:\O\to \O_0$
be the projection functor
(or its extension to the derived categories).
 We can identify $\O_0$ with the category
of modules over $End(P_e)$; the functor $\pi$ is then identified with
$X\mapsto Hom(P_e,X)$.

\medskip

{\bf Proposition.} The functor $\pi|_{\T}$ is fully faithful.

We will need the following standard fact:

\medskip

{\bf Lemma.} The socle of $\Delta_w$ and the cosocle of $\nabla_w$
are isomorphic to $L_e$.\ \ \

\medskip

{\it Proof} of Lemma. Let us prove the
statement about  $\Delta_w$; the one about
$\nabla_w$ then follows by Verdier duality. We
argue by induction in the length $\ell(w)$. If
$w=e$ there is nothing to prove, and if
$\ell(w)=1$ then the statement follows from
the existence of a non-split exact sequence
$$
0\to \Delta_e\to \Delta_w \to L_w\to 0 \tag 2.1.1
$$
of perverse sheaves on ${\Bbb P}^1$.

Assume now that $w=w's$, where $s$ is the simple reflection
corresponding to a simple root  $\alpha$,
 and $\ell(w)>\ell(w')$. Let $X^\alpha$
be the corresponding partial flag variety, and
 $pr_\alpha:X \to X^\alpha$ be the projection; thus $pr_\alpha$
is a fibration with projective lines as fibers. Set
$X^\alpha_w=pr_\alpha(X_w)$; $X _w'=pr_\alpha^{-1}
(X_w^\alpha)$, and let $i_w^\alpha: X_w^\alpha
\imbed X^\alpha$, $i_w':X_w'\imbed X$ be the
embeddings.
 Then $i_w'$ is an affine morphism because it
is a base change of the affine morphism
$i_w^\alpha$. Hence the functor $i_{w!}'$ is
exact. The fibration $pr_\alpha$ is trivial
over $pr(X_w)$, so we have $X_w'\cong \Pone
\times X_{w'}$.
 Applying the functor $i'_{w!}\circ pr_1^*
[\ell(w)-1]$ to (2.1.1) (where $pr_1:X_w'\to \Pone$ is the
projection)
we get an exact sequence in $\O$
$$0\to i_{w'!}(M_1) \to i_{w!} (M_2 ) \to 
i'_{w!}(M_3)\to 0,$$ where $M_1$, $M_2$, $M_3$
are constant perverse sheaves on the
corresponding varieties.  Let $L_u\subset
\Delta_w$ be a simple subobject. Suppose first
that the composition $L\to i'_{w!}(M_3)$ is
nonzero. It is easy to see that this only can
happen if $\ell(u\cdot s_\alpha)< \ell(u)$ so
that $L_u = pr_\alpha^*(L')[1]$ for a certain
irreducible perverse sheaf $L'$ on
$X^\alpha$. We arrive to a contradiction since
 $Hom(L_u,i_{w!}(M_2))=
Hom\left(L'[1], pr_{\alpha
\bullet}i_{w!}(M_2)\right) =$ $
Hom\left(L'[1],i^\alpha_{w!}(M_4[-1])\right)=
Ext^{-2}\left(L',
i^\alpha_{w!}(M_4)\right)=0$ (here $M_4$ is a
rank 1 perverse constant  sheaf on
$X_w^\alpha$; and we use the notation
 $ f_{\bullet}:= f_*=f_!$ for a proper morphism $f$).
  Thus we have 
$L\subset i_{w'!}(M_1 )$, so we
get the desired statement by induction. \hfill$\square$

\medskip

{\it Proof} of  Proposition.
 Let $\A$ be an abelian category, and $\B\subset\A$ be a Serre
subcategory. Define the left and right orthogonals to $\B$ in $\A$
by
$$^\perp \B= 
\{ A\in \A \ |\
Hom(A,X)=0 \ \ \forall X\in \B\}$$
$$ \B^\perp =
\{ A\in \A \ |\
Hom(X,A)=0 \ \ \forall X\in \B\}.$$ It follows
from the definitions that if  $A\in\, ^\perp
\B$ and $B\in \B^\perp$, then
$
Hom_\A(A,B)\iso Hom_{\A/\B}(A,B).$
The lemma implies that $\Delta_w\in
\O_{>0}^\perp$, $\nabla_w\in\, ^\perp\O_{>0}$
for all $w$. Hence $\T\subset \,
^\perp\O_{>0}\cap  \O_{>0}^\perp$, so our
proposition is proved. \hfill$\square$

\medskip

{\bf 2.2.} We recall the intertwining functors (Radon transforms)
acting on $D$. Let  $\ell(w)=\dim (X_w)$ be the length function.
For $w\in W$ let $X^2_w\subset X^2$ be the $G$-orbit corresponding
to $w$ (thus $X^2_w= G(X_e\times X_w)$). Let $pr_i^w:X^2_w\to X$
for $i=1,2$ be the projections. Set $R_w^?(X)=pr_{2?} pr_1^*(X)[\ell(w)]$,
where $?=!$ or $*$. We need a standard

{\bf Fact.} For $?=!$ or $*$ we have:

(a) $R_{w_1}^?\circ R_{w_2}^?\cong
R_{w_1w_2}^?$ provided $\ell(w_1w_2)=\ell(w_1)+\ell(w_2)$.

(b) $R_w^!\circ R_{w^{-1}}^*\cong id \cong
R_{w^{-1}}^* \circ R_w^!$.

(c) $\pi\circ R_w^? \cong \pi$.

\medskip

{\it Proof.} (a) and (b) are well known (see
e.g.~[BB]). 
 Using (a) we see that it is enough to check (c)
for $w$ of length 1; so assume that $w=s_\alpha$ is a simple reflection.
 We treat the case
$?=!$, the other case follows.
Let $\overline X_{s_\alpha}^2$ be the closure of $X_{s_\alpha}^2$, and let
$\overline{pr}_1$, $\overline{pr}_2:\overline X_{s_\alpha}^2\to X$ be the
projections. Thus $\overline{pr}_1$, $\overline{pr}_2$ are fibrations
with fiber $\Pone$. For $M\in D$ we have a canonical exact triangle
$$ \delta_*(M) \to i_!pr_1^* M [1]\to \overline{pr}_1^*M[1]$$
where $\delta:X\to X^2$ is the diagonal embedding, and $i:X_{s_\alpha}^2
\imbed X^2$ is the embedding. Applying 
$\overline{pr}_{2!}$ to it,
we see that it suffices to check that
 $$\pi( \overline{pr}_{2!}
\overline{pr}_1^*M)=0 \tag 2.2.1$$
 This is clear since 
$\overline{pr}_{2!} \overline{pr}_1^*M=
pr_\alpha^* pr_{\alpha !}M$ (we use notation of
the proof of Lemma 2.1). Indeed, the pull-back
functor from $X_\alpha$ identifies irreducible
perverse sheaves constant along the Schubert
stratification on $X_\alpha$ with irreducible
objects of $\CO$ constant along the fibers of
$p_\alpha$, so $L_e$ cannot occur in
$pr_\alpha^* pr_{\alpha !}M$.   
 \hfill$\square$

\medskip

{\bf 2.3.} The following result, inspired  by
W.~Soergel's article [S2], appears in [BG] (see
{\it loc.cit.}~Theorem 6.10(i)); it was also 
found independently by R.~Rouquier
(unpublished). We include a proof for the
reader's convenience.

Let $w_0\in W$ be the longest element.

{\bf Proposition.} 
We have $R_{w_0}^!(I_w)\cong T_{ww_0}$; $R_{w_0}^!(T_w)=P_{ww_0}$.

{\bf Lemma.} Assume we are in the situation of 
1.5; denote by $M_\nu$ a constant
perverse sheaf of rank 1 on $X_\nu$. Let $M\in
D(X,\{X_\nu\})$ be any perverse sheaf.

If $Ext^a (M, i_{\nu !}(M_\nu))=0$ for every
 $a>0$ and $\nu$, then $M$ is
projective.

If  $Ext^a ( i_{\nu *} (M_\nu), M)=0$ for every
$a>0$ and $\nu$, then $M$ is
injective.

{\it Proof} of  Lemma. We prove the first
statement, the second one is
 similar.
We will say that an object $A$
of a triangulated category is filtered by objects $B_i$
if there exist objects $A_0=0, \, A_1,\dots, A_n=A$ and exact triangles
$A_{i-1} \to A_i\to B_i$. Then the definition of the perverse
$t$ structure implies that any perverse sheaf $N$ is filtered (as
an object of the triangulated category $D(X,\{X_\nu \} )$)
by objects of the form $i_{\nu !}(M_\nu)[d]$,
$d\leq 0$. Thus the condition implies that
$Ext^a (M,N)=0$ for all $a>0$.
\hfill$\square$

{\it Proof} of  Proposition.
  We prove the first isomorphism, the second one is
 similar.
Let us first see that  $R_{w_0}^!(T_w)$ is a projective object of $\O$.
By Fact above we have $$R_{w_0}^! (\nabla _w)=R_{w_0w}^! \circ R_{w^{-1}}^!
(R_w^* (\Delta_e))=R_{w_0}^!(\Delta_e)=\Delta_{w_0w}.$$
It follows that $R_{w_0}^! (\nabla _w)\in D$ is filtered by the objects
$\nabla_w$, in particular, it lies in $\O$. We
also have
$$Ext^a (R_{w_0}^!(T_w), \Delta_w)=
Ext^a (R_{w_0}^!(T_w), R_{w_0}^! (\nabla_
{w_0w}))=Ext^a (T_w, \nabla_{w_0w})=0
$$
 for $ a>0$,
where the last equality follows from the fact that $T_w$ is filtered by
objects $\Delta_u$, and $Ext^{>0}(\Delta_u,
\nabla_v)=0$. Thus, by Lemma, we see that 
$R_{w_0}^!(T_w)$ is projective. Moreover,
 $R_{w_0}^!(T_w)$ is indecomposable, and it
follows from the above that it is filtered by
objects of the form $\Delta_{w_i}$ where
$w_1=w_0w$, and
$w_i\succ w_0w$ for $i>1$. It follows that  $R_{w_0}^!(T_w)\cong P_{w_0w}$.
We have proved the second isomorphism; by Verdier duality it
implies that $R_{w_0}^*(T_w)\cong I_{w_0 w}$;
and applying $R_{w_0}^!$ to both sides we get
the first isomorphism.
\hfill$\square$

\medskip

{\bf 2.4 Corollary.} (Soergel's Struktursatz,
[S1], p.433) The functors $\pi|_{\I}$,
$\pi|_{\P}$ are fully faithful.
\hfill$\square$

\medskip

{\bf 2.5.} Recall some definitions of Bondal
and Kapranov [BK]. Let $\D$ be a $k$-linear
category such that
$Hom(X,Y)$ is a finite-dimensional vector
space for all $X,Y\in \D$. Suppose that $\D$
admits an endofunctor $S$ equipped with
a natural isomorphism $\alpha : Hom
(X,S(Y))\iso Hom (Y,X)^*$. Such $(S,\alpha )$
is evidently unique.\footnote{And for given $S$
all possible $\alpha$ form a torsor with
respect to an obvious action of the group of
automorphisms of the identity functor
$Id_\D$.} It is called  {\it the Serre
functor} if
$S$ is actually  an auto-equivalence of
$\D$.\footnote{Additional requirements on $S$
imposed in the definition  of the Serre functor
in [BK] are actually redundant,
see the proof of Proposition 3.4  in {\it
loc.~cit.}} If $\D$ is a triangulated
category, then $S$ is naturally a triangulated
functor.

Let us return to our situation. The Serre
functor on
$D$ exists by the results [BK] (compare [BK],
Corollary 3.5 with either Theorem 2.11 or
Corollary 2.10 in {\it loc.~cit.}).  In fact,
 the bounded derived category
$D^b (A)$ has the Serre functor whenever $A$ is
an Artinian abelian category of
finite homological dimension having enough
projectives and finitely many isomorphism
classes of irreducible objects.

 The following result was conjectured by
Kapranov:

{\bf Proposition.} The Serre functor $S$ for
$D$ is isomorphic to $(R_{w_0}^*)^2$  as a
triangulated functor.

{\it Proof.} It takes two steps: 

(i) Our functors send $\P$
to $\I$ and their restrictions to
$\P$ are isomorphic.

(ii) Any isomorphism of
functors $(R^*_{w_0})^2 |_\P \iso S|_\P$
extends in a canonical way to an isomorphism of
triangulated functors $(R^*_{w_0})^2  \iso
S$.

{\it Proof of (i).} 
Notice that for
each $w\in W$ one has $S(P_w )\cong I_w \cong
(R_{w_0}^*)^2(P_w)$ as follows, respectively,
from [BK] 3.2(3) and  2.3.

We will prove that any isomorphism $S(P_e
)\iso  (R_{w_0}^*)^2 (P_e )$
extends uniquely to an isomorphism of functors
$S|_\P \iso (R_{w_0}^*)^2 |_\P$. 

According to 2.4, we can replace our functors
by their composition with $\pi$.  By
Fact 2.2(c), one has
$\pi (R_{w_0}^*)^2 \cong \pi$, so we can
reformulate our claim as follows: Any
isomorphism $\alpha : \pi S(P_e )\iso \pi (P_e
)$ extends uniquely to an isomorphism of
functors $\pi S|_\P \iso \pi|_\P$. Since $\pi
= Hom (P_e ,\cdot )$, it suffices to check that
 $\alpha$ commutes with the
action of $End(P_e )$. As follows
from the definition of the Serre functor, $S$
commutes with any endomorphism of the identity
functor  $Id_D$. Now any endomorphism of $P_e$
comes from an endomorphism of $Id_D$, as
follows from  2.4 and
commutativity of $End(P_e )$ established in 
[S1], Lemma 5
on p.~430,\footnote{In fact, Soergel's
``Endomorphismensatz''  ({\it
loc.~cit.}~p.~428; see also [B]) provides a
very explicit description of
$End(P_e)$, see Remark (ii)
below.} and we are done.

{\it Proof of (ii).}
For a $k$-linear functor $\phi : \P \to\I$
let $C(\phi ) : C^b (\P )\to C^b (\I
)$ be its DG extension to the category of
bounded complexes and $D(\phi )$ the
triangulated endofunctor of $D =K^b (\P )=K^b
(\I )$ defined by $C(\phi )$. We have seen
that the restrictions of $S$ and
$(R^*_{w_0})^2$ to
$\P$ are isomorphic functors $\CP \to \CI$.
We will show that there are {\it canonical}
identifications of triangulated endofunctors $S
\iso D(S|_\P )$,
$(R^*_{w_0})^2 \iso D((R^*_{w_0})^2 |_\P )$;
this yields (ii).

The statement about $S$ is clear. Indeed,
since $S$ is the Serre functor, we have a
natural isomorphism
$Hom (X,S(Y))\iso Hom (Y,X)^*$ for $X\in \O$,
$Y\in\P$. It extends canonically to a natural
DG isomorphism $Hom (X,S(Y))\iso Hom (Y,X)^*$
for $X\in C^b (\O )$,
$Y\in C^b (\P )$, which makes $D(S|_\P )$ the
Serre functor. 

Consider  $(R^*_{w_0})^2$. This is the
restriction to $D$ of the endofunctor
$(R^*_{w_0})^2 $ of the derived category $D
(X)$ of constructible complexes on $X$. The
latter functor has ``geometric origin" hence
it lifts canonically to a triangulated 
endofunctor $(R^{*F}_{w_0})^2$ of the filtered
derived category $D F (X)$ of finitely
filtered constructible complexes. Recall (see
[BBD] 3.1) that there is a canonical fully
faithful embedding $C^b (\CM (X)) \hra DF (X)$
whose essential image consists of those
filtered complexes $P$ that $\gr^i P \in \CM
(X)[-i]$ for any $i$ (the inverse
functor identifies such $P$ with the complex of
perverse sheaves $..\to \gr^i P \to
\gr^{i+1} P\to ..$ where the differential is
the third side of the triangle $\gr^{i+1} P \to
P_{i}/P_{i+2}
\to\gr_{i}P$). The  equivalence
$D^b (\O )\iso D\subset D (X)$ from 1.5.1
comes from the composition $C^b (\O )\hra C^b
(\CM (X))\hra D F (X) \to D (X)$ where 
the third arrow is the forgetting of
filtration functor. Now the restriction of
$(R^{*F}_{w_0})^2$ to $C^b (\P )\subset C^b
(\CM (X))$ is  $C^b (\P ) \to
C^b (\I ) \subset C^b (\CM (X))\subset D F
(X)$ where the arrow is $C(
(R^*_{w_0})^2 |_\P )$. Since
$(R^{*F}_{w_0})^2$ lifts $(R^{*}_{w_0})^2$, we
are done.  \hfill$\square$

\medskip

{\bf 2.6 Remarks.} (i) Let $\tilde W$ be 
 the braid group  associated
 to the root system of $G$, and for $w\in W$ let $\tilde w\in \tilde W$
be its canonical (minimal length) lifting. 
According to 2.2, for $?=*,!$
the map
$w\mapsto  R_w^?$ extends to a weak action of
$\tilde W$ on $D$  (extending it to a strong
action in the sense of [D] requires more work;
this is done in [R]).  Notice that
$\tilde{w_0}^2$ is a central element in
$\tilde W$. This conforms with the general
fact that for any triangulated category
$\D$ with a Serre functor $S$, and any other functor $F: \D\to \D$
which admits a left adjoint $LF$ we have a canonical isomorphism
$$F\circ S \cong S\circ L(LF)$$
(where $L(LF)$ denotes the left adjoint to $LF$). In particular,
if $F$ is invertible we have $LF\cong F^{-1}$, so $L(LF)\cong F$,
i.e. $F$ commutes with $S$.

 (ii) In the step (i) of the proof of
Proposition we have shown that the set of
isomorphism of functors
$(R_{w_0}^*)^2 |_\P \iso S|_\P$
identifies canonically with the
$Z^\times$-torsor of invertible elements in
the $Z$-module $K:= Hom ((R_{w_0}^*)^2 (P_e
),  S(P_e ))$ where $Z:= End (P_e )$. Now
$K\iso  Hom (P_e , (R_{w_0}^!)^2 (P_e ))^* \iso
Hom (\pi (P_e ), \pi((R_{w_0}^!)^2 (P_e
)))^*$ which  equals $Hom (\pi (P_e ),\pi
(P_e ))^* \iso Hom (P_e , P_e )^* \iso Z^*$ 
(the $k$-linear dual to $Z$) by  2.2(c). Thus
we have a canonical isomorphism of
$Z$-modules
$K\iso Z^*$.
Accordnig to [S1], [B],
 there is a canonical isomorphism of
algebras
$Z\iso H^* (X^\vee )$, where $X^\vee$ is the
flag space for the Langlands dual group
$G^\vee$.\footnote{Actually [S1], [B] work with
modules over the enveloping algebra, so one
has to invoke the localization theorem to
derive the computation of $End(P_e)$ from their
results. There is an equivalent, purely
topological, construction (see [BGS], p.~525)
of the morphism
$H^\cdot (X^\vee )\to A$. One knows that
$H^\cdot (X^\vee )$ is generated by $H^2
(X^\vee )$, and the Chern class for the
$T^\vee$-torsor
$G^\vee /N^\vee$ over
$X^\vee$ provides a canonical identification
$\frak t  \iso H^2 (X^\vee
)$, where
$\frak t$ is the Cartan algebra of $G$. So our
morphism is determined by a linear map
$\frak t \to Z$.  Our perverse sheaves
are monodromic (of unipotent monodromy) with
respect to the action of (any) maximal torus
$T\subset G$ on
$X$. Now $\frak t \to Z$ is  the logarithm
of the monodromy map.}  So the trace map
$H^* (X^\vee )\to k$  provides a canonical
generator of the
$Z$-module $Z^*$. It yields  a canonical
isomorphism of functors $(R_{w_0}^*)^2 |_\P
\iso S|_\P$ hence, by step (ii) of the proof of
Proposition, an identification of the 
triangulated functors
$(R_{w_0}^*)^2 \iso S$.

\medskip

{\bf Acknowledgment.}
We thank Michael Finkelberg who taught us the
theory of tilting objects.
 We are grateful to  
Raphael Rouquier and Victor Ginzburg
for help with  references; we are also much obliged
to Rouquier for pointing out a mistake
 in the first 
version of the text. R.B. thanks Independent Moscow University where
part of this work was done.

\medskip

{\bf References.}

[B1] A.~Beilinson, How to glue perverse sheaves,
Springer LN 1289, 42--51.

[B2] A.~Beilinson, On the derived category of perverse
sheaves, ibid., 27--41.

[BB] A.~Beilinson, J.~Bernstein,
 A generalization of
Casselman's submodule theorem, in:
 ``Representation theory of reductive groups'' (Park City,
Utah, 1982), 35--52, Progr. Math., 40, Birkh\" auser Boston, Boston, Mass., 
1983.

[BBD] A.~Beilinson, J.~Bernstein, P.~Deligne,
Faisceaux pervers, Ast\'erisque 100, 1982.

[BG] A.~Beilinson, V.~Ginzburg, Wall-crossing functors and D-modules,
Repr. Theory 3(1999), 1-31.

[BGS] A.~Beilinson, V.~Ginzburg, W.~Soergel,
 Koszul duality patterns in
representation theory, J.~Amer.~Math.~Soc.~9,
no.~2, 473--527, 1996.

[B] J.~Bernstein, 
Trace in categories; in 
``Operator algebras, unitary representations,
 enveloping algebras, and invariant theory'' (Paris, 1989), 417--423, 
Progr. Math., 92, 
Birkh\" auser Boston, Boston, MA, 1990. 

[BK] A.~Bondal, M.~Kapranov, Representable
functors, Serre functors, and mutations.
Izv.~Acad.~Nauk SSSR. Ser.~Mat., 53,
1183--1205, 1989; English
translation: Math.~USSR.~Izv., 35, 519--541,
1990.

[D] P.~Deligne, 
Action du groupe des tresses sur une cat\' egorie, 
Invent. Math. 128 (1997), no. 1, 159--175.

[Di] C.~Dickens, {\it Oliver Twist;} chapter 18. 

[M] O.~Mathieu, Tilting modules and their applications.
Analysis on homogeneous spaces and representation theory
of Lie groups, Okayama--Kyoto (1997), 145--212, Adv. Stud.
Pure Math., 26, Math. Soc. Japan, Tokyo, 2000.

[MV] R.~MacPherson, K.~Vilonen, Elementary construction of
perverse sheaves, Invent.~Math.~84, 403--435, 1986.

[Ri]
C.M. Ringel, 
The category of modules with good filtrations over a quasi-hereditary 
algebra has almost split sequences.
Math.~Z.~208 (1991), no. 2, 209--223.

[R] R.~Rouquier, Action du groupe de tresses
sur la cat\' egorie d\' eriv\' ee de la vari\'
et\' e de drapeaux,
 preprint, available at:
http://www.math.jussieu.fr/$\sim$rouquier/preprints/flag.dvi

[S1] W.~Soergel, Kategorie $\O$, perverse
Garben und Moduln \"uber den Koinvarianten zur
     Weylgruppe, J.~Amer.~Math.~Soc.~3,
no.~2, 421--445, 1990.

[S2] W.~Soergel,
Character formulas for tilting modules over 
Kac-Moody algebras, 
Represent.~Theory~2 (1998), 432--448.
\end


**********************************

recall the convolution functor
$D\times D_B(X)\to D$, $(X,Z) \mapsto X*Z$
 where $D_B$ denotes the equivariant derived category
of sheaves (or $D$-modules) on $X$. The functor
 $R_w^!$ is identified with the functor $X\mapsto X*  \Delta _w$
(where we use the same notation $\Delta_w$ for objects of
two different categories: $D$ and $D_B(X)$). For $w$ of length 1
we have an exact sequence $0\to L_e \to \Delta_w \to L_w$, which
yields a functorial exact triangle $X\to R_w^!(X) \to X*L_w$.
However, it is easy to see that $\pi(X*L_w)=0$ for all $X$ provided $w\ne e$.
This shows the statement concerning $R_w^!$, and the one concerning
$R_w^*$ is proved dually.

***********************************

 In this remark we address
the natural question of whether one can choose
a canonical isomorphism $(R_{w_0}^*)^2\cong S$.
In fact, we claim that there a canonical
isomorphism
$$ A^*
\mathop\otimes\limits_Z (R^*_{w_0})^2 \cong S \tag 2.3.1
$$ where $Z=End(P_e)$ is the center of $D$, and 
$Z^*$ is the linear dual to $Z$ (the cofree module); since $Z
=End(P_e)\cong 
H^*(G/B)$ (where the last isomorphism is due to [S1], and the first
one (also proved in [S1]) follows from Corollary 2.2)
is a Frobenius algebra, the left hand side in 2.3.1 is indeed isomorphic
to $(R^*_{w_0})^2 $.

To define 2.3.1
notice that
Proposition 2.2 (or Corollary 2.2) imply that to specify an isomorphism
2.3.1 it is enough to specify  the induced isomorphism 
$$ Z^*
\mathop\otimes\limits_Z
(R_{w_0})^2(P_e)\cong  S (P_e) \tag 2.3.2 $$ The proof  of Fact 2.2(c) 
shows that there exists a canonical isomorphism $\pi \circ R_w\cong \pi$,
and hence a canonical isomorphism $R_w(P_e)\cong P_e$, which also 
yields a canonical isomorphism $R_{w_0}^2(P_e)\cong P_e$. Thus to fix $(2.3.2)$
we need to fix an isomorphism $S(P_e)\cong   Z^*
\mathop\otimes\limits_Z
P_e$. Finally, such an
isomorphism is uniquely defined by the induced isomorphism
$ Hom (P_e, S(P_e))\cong Z^*\mathop\otimes\limits_Z
Hom(P_e,P_e)=Z^*\mathop\otimes\limits_Z Z=Z^*$. Here the right hand side 
is identified with the dual space to $Z=End(P_e)$ by the 
definition of a Serre functor.


 to $End(P_e)$ is an isomorphism. It also shows
that every autoequivalence of $\P$ which (i)
sends $P_w$ to
$P_w$,  and (ii) acts trivially on $Z$,
is isomorphic to identity. Thus it suffices
 to see that $(R_{w_0}^!)^2\circ S$ satisfies
(i), (ii).

To see  property (ii) we check that both $S$ and $R_{w}^!$
act trivially on the center of $D$. 
Here the claim about $S$ is a general simple
property of the Serre functor. The claim
about
$R_w^!$ follows from Fact (c) in 2.2, and the
fact that $Z\iso End(P_e)$.